\documentclass[12pt,twoside]{article}
\usepackage{amsmath,amssymb}
\usepackage{amsthm}
\usepackage{upref}
\numberwithin{equation}{section}

\newtheorem*{theorem*}{Theorem}
\theoremstyle{plain}
\newtheorem{theorem}{Theorem}[section]
\newtheorem{proposition}[theorem]{Proposition}
\newtheorem{definition}[theorem]{Definition}

\theoremstyle{definition}
\newtheorem{remark}[theorem]{Remark}
\newtheorem{problem}[theorem]{Problem}

\begin{document}

\newcommand{\eq}{equation}
\newcommand{\real}{\ensuremath{\mathbb R}}
\newcommand{\comp}{\ensuremath{\mathbb C}}
\newcommand{\rn}{\ensuremath{{\mathbb R}^n}}
\newcommand{\tn}{\ensuremath{{\mathbb T}^n}}
\newcommand{\rnp}{\ensuremath{\real^n_+}}
\newcommand{\rnn}{\ensuremath{\real^n_-}}
\newcommand{\Rn}{\ensuremath{{\mathbb R}^{n-1}}}
\newcommand{\Zn}{\ensuremath{{\mathbb Z}^{n-1}}}
\newcommand{\no}{\ensuremath{\nat_0}}
\newcommand{\ganz}{\ensuremath{\mathbb Z}}
\newcommand{\zn}{\ensuremath{{\mathbb Z}^n}}
\newcommand{\zom}{\ensuremath{{\mathbb Z}_{\Om}}}
\newcommand{\zOm}{\ensuremath{{\mathbb Z}^{\Om}}}
\newcommand{\As}{\ensuremath{A^s_{p,q}}}
\newcommand{\Bs}{\ensuremath{B^s_{p,q}}}
\newcommand{\Fs}{\ensuremath{F^s_{p,q}}}
\newcommand{\Fsr}{\ensuremath{F^{s,\rloc}_{p,q}}}
\newcommand{\nat}{\ensuremath{\mathbb N}}
\newcommand{\Om}{\ensuremath{\Omega}}
\newcommand{\di}{\ensuremath{{\mathrm d}}}
\newcommand{\sn}{\ensuremath{{\mathbb S}^{n-1}}}
\newcommand{\Ac}{\ensuremath{\mathcal A}}
\newcommand{\Acs}{\ensuremath{\Ac^s_{p,q}}}
\newcommand{\Bc}{\ensuremath{\mathcal B}}
\newcommand{\Cc}{\ensuremath{\mathcal C}}
\newcommand{\cc}{{\scriptsize $\Cc$}${}^s (\rn)$}
\newcommand{\ccd}{{\scriptsize $\Cc$}${}^s (\rn, \delta)$}
\newcommand{\Fc}{\ensuremath{\mathcal F}}
\newcommand{\Lc}{\ensuremath{\mathcal L}}
\newcommand{\Mc}{\ensuremath{\mathcal M}}
\newcommand{\Ec}{\ensuremath{\mathcal E}}
\newcommand{\Pc}{\ensuremath{\mathcal P}}
\newcommand{\Efr}{\ensuremath{\mathfrak E}}
\newcommand{\Mfr}{\ensuremath{\mathfrak M}}
\newcommand{\Mbf}{\ensuremath{\mathbf M}}
\newcommand{\Dbb}{\ensuremath{\mathbb D}}
\newcommand{\Lbb}{\ensuremath{\mathbb L}}
\newcommand{\Pbb}{\ensuremath{\mathbb P}}
\newcommand{\Qbb}{\ensuremath{\mathbb Q}}
\newcommand{\Rbb}{\ensuremath{\mathbb R}}
\newcommand{\vp}{\ensuremath{\varphi}}
\newcommand{\hra}{\ensuremath{\hookrightarrow}}
\newcommand{\supp}{\ensuremath{\mathrm{supp \,}}}
\newcommand{\ssupp}{\ensuremath{\mathrm{sing \ supp\,}}}
\newcommand{\dist}{\ensuremath{\mathrm{dist \,}}}
\newcommand{\unif}{\ensuremath{\mathrm{unif}}}
\newcommand{\ve}{\ensuremath{\varepsilon}}
\newcommand{\vk}{\ensuremath{\varkappa}}
\newcommand{\vr}{\ensuremath{\varrho}}
\newcommand{\pa}{\ensuremath{\partial}}
\newcommand{\oa}{\ensuremath{\overline{a}}}
\newcommand{\ob}{\ensuremath{\overline{b}}}
\newcommand{\of}{\ensuremath{\overline{f}}}
\newcommand{\LA}{\ensuremath{L^r\!\As}}
\newcommand{\LcA}{\ensuremath{\Lc^{r}\!A^s_{p,q}}}
\newcommand{\LcdA}{\ensuremath{\Lc^{r}\!A^{s+d}_{p,q}}}
\newcommand{\LcB}{\ensuremath{\Lc^{r}\!B^s_{p,q}}}
\newcommand{\LcF}{\ensuremath{\Lc^{r}\!F^s_{p,q}}}
\newcommand{\Lf}{\ensuremath{L^r\!f^s_{p,q}}}
\newcommand{\La}{\ensuremath{\Lambda}}
\newcommand{\Lob}{\ensuremath{L^r \ob{}^s_{p,q}}}
\newcommand{\Lof}{\ensuremath{L^r \of{}^s_{p,q}}}
\newcommand{\Loa}{\ensuremath{L^r\, \oa{}^s_{p,q}}}
\newcommand{\Lcoa}{\ensuremath{\Lc^{r}\oa{}^s_{p,q}}}
\newcommand{\Lcob}{\ensuremath{\Lc^{r}\ob{}^s_{p,q}}}
\newcommand{\Lcof}{\ensuremath{\Lc^{r}\of{}^s_{p,q}}}
\newcommand{\Lca}{\ensuremath{\Lc^{r}\!a^s_{p,q}}}
\newcommand{\Lcb}{\ensuremath{\Lc^{r}\!b^s_{p,q}}}
\newcommand{\Lcf}{\ensuremath{\Lc^{r}\!f^s_{p,q}}}
\newcommand{\id}{\ensuremath{\mathrm{id}}}
\newcommand{\tr}{\ensuremath{\mathrm{tr\,}}}
\newcommand{\trd}{\ensuremath{\mathrm{tr}_d}}
\newcommand{\trL}{\ensuremath{\mathrm{tr}_L}}
\newcommand{\ext}{\ensuremath{\mathrm{ext}}}
\newcommand{\re}{\ensuremath{\mathrm{re\,}}}
\newcommand{\Rea}{\ensuremath{\mathrm{Re\,}}}
\newcommand{\Ima}{\ensuremath{\mathrm{Im\,}}}
\newcommand{\loc}{\ensuremath{\mathrm{loc}}}
\newcommand{\rloc}{\ensuremath{\mathrm{rloc}}}
\newcommand{\osc}{\ensuremath{\mathrm{osc}}}
\newcommand{\pr}{\pageref}
\newcommand{\wh}{\ensuremath{\widehat}}
\newcommand{\wt}{\ensuremath{\widetilde}}
\newcommand{\ol}{\ensuremath{\overline}}
\newcommand{\os}{\ensuremath{\overset}}
\newcommand{\Li}{\ensuremath{\overset{\circ}{L}}}
\newcommand{\Ai}{\ensuremath{\os{\, \circ}{A}}}
\newcommand{\Ci}{\ensuremath{\os{\circ}{\Cc}}}
\newcommand{\Si}{\ensuremath{\os{\, \circ}{S}}}
\newcommand{\dom}{\ensuremath{\mathrm{dom \,}}}
\newcommand{\SA}{\ensuremath{S^r_{p,q} A}}
\newcommand{\SB}{\ensuremath{S^r_{p,q} B}}
\newcommand{\SF}{\ensuremath{S^r_{p,q} F}}
\newcommand{\Hc}{\ensuremath{\mathcal H}}
\newcommand{\Nc}{\ensuremath{\mathcal N}}
\newcommand{\Lci}{\ensuremath{\overset{\circ}{\Lc}}}
\newcommand{\bmo}{\ensuremath{\mathrm{bmo}}}
\newcommand{\BMO}{\ensuremath{\mathrm{BMO}}}
\newcommand{\cm}{\\[0.1cm]}
\newcommand{\Aa}{\ensuremath{\os{\, \ast}{A}}}
\newcommand{\Ba}{\ensuremath{\os{\, \ast}{B}}}
\newcommand{\Fa}{\ensuremath{\os{\, \ast}{F}}}
\newcommand{\Aas}{\ensuremath{\Aa{}^s_{p,q}}}
\newcommand{\Bas}{\ensuremath{\Ba{}^s_{p,q}}}
\newcommand{\Fas}{\ensuremath{\Fa{}^s_{p,q}}}
\newcommand{\Ca}{\ensuremath{\os{\, \ast}{{\mathcal C}}}}
\newcommand{\Cas}{\ensuremath{\Ca{}^s}}
\newcommand{\Car}{\ensuremath{\Ca{}^r}}
\newcommand{\bl}{$\blacksquare$}

\begin{center}
{\Large Mapping properties of the Fourier transform in spaces with dominating mixed smoothness}
\\[1cm]
{Hans Triebel}
\\[0.2cm]
Institut f\"{u}r Mathematik\\
Friedrich--Schiller--Universit\"{a}t Jena\\
07737 Jena, Germany
\\[0.1cm]
email: hans.triebel@uni-jena.de
\end{center}

\begin{abstract}
This paper deals  with continuous and compact mappings of the Fourier transform in function spaces with dominating mixed smoothness.
\cm
{\bfseries Keywords:} Fourier transform, function spaces, dominating mixed smoothness
\cm
{\bfseries 2020 MSC:} 46E35
\end{abstract}

\section{Introduction and main assertions}   \label{S1} 
Let $\nat$ be the set of all natural numbers and $\no = \nat \cup \{0 \}$. Let $\rn$ be the Euclidean $n$-space, where  
$n \in \nat$. Put $\real = \real^1$, whereas $\comp$ is the complex plane. Let $L_p (\rn)$ with $0<p \le \infty$, be the standard quasi-Banach space with respect to the Lebesgue measure in \rn, quasi-normed by
\begin{\eq}   \label{1.1}
\| f \, | L_p (\rn) \| = \Big( \int_{\rn} | f(x)|^p \, \di x \Big)^{1/p}
\end{\eq}
with the natural modification if $p=\infty$. As usual, $\ganz$ is the set of all integers and $\zn$  where $n \in \nat$, denotes the lattice of
all points $m = (m_1, \ldots, m_n )\in \rn$ with $m_j \in \ganz$. Let $\nat^n_0$, where $n \in \nat$, be the set of all multi-indices,
\begin{\eq}   \label{1.2}
\alpha = (\alpha_1, \ldots, \alpha_n ) \quad \text{with $\alpha_j \in \no$ and $|\alpha| = \sum^n_{j=1} \alpha_j$.}
\end{\eq}
Let as usual $\pa_j =\pa/\pa x_j$, $\pa^m_j = \pa^m/\pa x^m_j$, $m\in \no$, ($\pa^0_j f =f$) and $D^\alpha = \pa_1^{\alpha_1} \cdots \pa^{\alpha_n}_n$, $\alpha \in \nat^n_0$.

Let $1<p<\infty$ and $r\in \no$. Then $S^r_p W(\rn)$
are the classical Sobolev spaces with dominating mixed smoothness which can be equivalently normed by
\begin{\eq}   \label{1.3}
\| f \, | S^r_p W (\rn) \| = \sum_{\substack{\alpha \in \nat^n_0, \\ 0\le \alpha_j \le r}} \| D^\alpha f \, | L_p (\rn) \| 
\sim \sum_{\substack{\alpha \in \nat^n_0, \\ \alpha_j \in \{0,r \}}} \| D^\alpha f \, | L_p (\rn) \|.
\end{\eq}
They were introduced by S.M. Nikol'skij in \cite{Nik62, Nik63a, Nik63b}. One may also consult the relevant parts of \cite{Nik77} (first edition 1969) and \cite{BIN75}. These spaces, being the dominating mixed counterpart of the classical isotropic Sobolev spaces $W^r_p (\rn)$, $r \in \no$, $1<p<\infty$, have been generalized and modified in many directions. The systematic Fourier--analytical 
approach of several versions of spaces with dominating mixed smoothness of type $\SA (\rn)$, $A\in \{B,F \}$, $r\in \real$ and $0<p,q
\le \infty$ ($p<\infty$ for $F$--spaces), being the dominating mixed counterpart of the related isotropic spaces $A^r_{p,q} (\rn)$, 
goes back to H.--J. Schmeisser, \cite{Schm80, Schm82} (his habilitation) and can also be found in \cite[Chapter 2]{ST87}. One may also
consult \cite[Section 2.1, pp.~80/81]{ST87} for further historical comments and references of the early history of these spaces. The 
next decisive step goes back to J. Vyb\'iral, \cite{Vyb06} (his PhD--theses) including atomic and wavelet representations of spaces 
with dominating mixed smoothness. We returned in \cite{T10} to this topic dealing with Haar bases and Faber bases in spaces with
dominating mixed smoothness applied to numerical integration and discrepancy (number--theoretical assertions about the distribution 
of points, especially in cubes in \rn). This has been complemented  in \cite{T19}, covering pointwise multipliers and corresponding
spaces in arbitrary domains in \rn. What follows might be considered as a further step in the theory of spaces with dominating mixed
smoothness, concentrating on mapping properties of the Fourier transform in suitable spaces. This is the dominating mixed counterpart
of related assertions for isotropic spaces obtained quite recently in \cite{Tri22, HST23} and the forthcoming book \cite{HST24} which
we describe now briefly.

Let $B^s_p (\rn) = B^s_{p,p} (\rn)$ with $s\in \real$ and $0< p,q \le \infty$ be the well--known special isotropic Besov spaces 
(recalled below in Definition \ref{D2.1} and Remark \ref{R2.2}). Let $S(\rn)$ be the usual Schwartz space  and $S'(\rn)$ be its dual,
the space of tempered distributions in \rn. If $\vp \in S(\rn)$ then
\begin{\eq}   \label{1.4}
\wh{\vp} (\xi) = (F \vp) (\xi) = (2 \pi)^{-n/2}  \int_{\rn} e^{-ix \xi} \vp(x) \, \di x, \quad \xi \in \rn,
\end{\eq}
denotes the the Fourier transform of \vp, whereas $F^{-1} \vp$ and  $\vp^\vee$ stand for the inverse Fourier transform given by the
right--hand side of \eqref{1.4} with $i$ in place of $-i$. Here $x \xi$ denotes the scalar product in \rn. Both $F$ and $F^{-1}$ are
extended to $S'(\rn)$ in the standard way. The mapping properties of the Fourier transform in the isotropic spaces $B^s_p (\rn)$ which
we have in mind can be summarized as follows. 
\cm
(i) {\em Let $1\le p \le 2$. Then there is a continuous  mapping}
\begin{\eq}   \label{1.5}
F: \quad B^{2n (\frac{1}{p} - \frac{1}{2}) + s_1}_p (\rn) \hra B^{-s_2}_p (\rn)
\end{\eq}
{\em if, and only if, both $s_1 \ge 0$ and $s_2 \ge 0$. This mapping is compact if, and only if, both $s_1 >0$ and $s_2 >0$.}
\cm
(ii) {\em Let $2 \le p \le  \infty$. Then there is a continuous mapping}
\begin{\eq}   \label{1.6}
F: \quad B^{s_1}_p (\rn) \hra B^{2n(\frac{1}{p} - \frac{1}{2}) -s_2}_p (\rn)
\end{\eq}
{\em if, and only if, both $s_1 \ge 0$ and $s_2 \ge 0$. This mapping is compact if, and only if, both $s_1 >0$ and $s_2 >0$.}
\cm
This coincides with \cite[Theorem 2.12]{HST24} based on \cite{HST23}, improving preceding assertions in \cite{Tri22}. It is the main
aim of the present paper to prove the following counterpart of these assertions for the special Besov spaces $S^r_p B(\rn)$, $1 \le
p \le \infty$, $r\in \real$, with dominating mixed smoothness as defined below in Definition \ref{D2.3} and Remark \ref{R2.4}.

\begin{theorem*} Let $n\in \nat$. Let $r_1 \in \real$ and $r_2 \in \real$.
\cm
{\em (i)} Let $1 \le p \le 2$. Then there is a continuous mapping
\begin{\eq}  \label{1.7}
F: \quad S^{\frac{2}{p}-1 +r_1}_p B (\rn) \hra S^{-r_2}_p B (\rn)
\end{\eq}
if, and only if, both $r_1 \ge 0$ and $r_2 \ge 0$. This mapping is compact if, and only if, both $r_1 >0$ and $r_2 >0$.
\cm
{\em (ii)} Let $2 \le p \le \infty$. Then there is a continuous mapping
\begin{\eq}  \label{1.8}
F: \quad S^{r_1}_p B(\rn) \hra S^{\frac{2}{p} -1 -r_2}_p  B (\rn)
\end{\eq}
if, and only if, both $r_1 \ge 0$ and $r_2 \ge 0$. This mapping is compact if, and only if, both $r_1 >0$ and $r_2 >0$.
\end{theorem*}

In Section \ref{S2} we collect further definitions and some prerequisites. The proof of the above Theorem is shifted to Section 
\ref{S3}.

\section{Definitions and prerequisites}   \label{S2}
Some basic notation were already introduced in Section \ref{S1}. This applies in particular to the Fourier transform and its inverse 
in the space of tempered distributions $S'(\rn)$ in the Euclidean $n$--space \rn. We fix our use of $\sim$ (equivalence) as already
used in \eqref{1.3}. Let $I$ be an arbitrary index set. Then
\begin{\eq}    \label{2.1}
a_i \sim b_i \quad \text{for} \quad i\in I \quad \text{(equivalence)}
\end{\eq}
for two sets of positive numbers $\{ a_i: i \in I \}$ and $\{ b_i: i\in I \}$ means that there are two positive numbers $c_1$ and $c_2$ such
that 
\begin{\eq}   \label{2.2}
c_1 a_i \le b_i \le c_2 a_i \qquad \text{for all} \quad i\in I.
\end{\eq}

Next we recall the definition of the isotropic function spaces $\As (\rn)$, $A\in \{B,F \}$, in $\rn$ and their dominating mixed 
counterparts $\SA (\rn)$.

Let $\vp_0 \in S(\rn)$ with
\begin{\eq}   \label{2.3}
\vp_0 (x) =1 \ \text{if $|x| \le 1$} \quad \text{and} \quad \vp_0 (x) =0 \ \text{if $|x| \ge 3/2$,}
\end{\eq}
and let
\begin{\eq}   \label{2.4}
\vp_k (x) = \vp_0 \big( 2^{-k} x\big) - \vp_0 \big( 2^{-k+1} x \big), \qquad x \in \rn, \quad k \in \nat.
\end{\eq}
Since
\begin{\eq}   \label{2.5}
\sum^\infty_{j=0} \vp_j (x) =1 \qquad \text{for $x \in \rn$},
\end{\eq}
the $\vp_j$ form a dyadic resolution of unity. The entire analytic functions  $(\vp_j \wh{f} )^\vee (x)$ make sense pointwise in $\rn$ for any $f\in S' (\rn)$.

\begin{definition}   \label{D2.1}
Let $\vp = \{ \vp_j \}^\infty_{j=0}$ be the above dyadic resolution of unity. 
\cm
{\em (i)} Let
\begin{\eq}   \label{2.6}
0<p \le \infty, \quad 0<q \le \infty, \quad s \in \real.
\end{\eq}
Then $\Bs (\rn)$ is the collection of all $f \in S' (\rn)$ such that
\begin{\eq}   \label{2.7}
\| f \, | \Bs (\rn) \|_{\vp} = \Big( \sum^\infty_{j=0} 2^{jsq} \| (\vp_j \widehat{f})^\vee | L_p (\rn) \|^q \Big)^{1/q} < \infty
\end{\eq}
$($with the usual modification if $q = \infty)$.
\cm
{\em (ii)} Let
\begin{\eq}   \label{2.8}
0<p<\infty, \quad 0<q\le \infty, \quad s \in \real.
\end{\eq}
Then $\Fs (\rn)$ is the collection of all $f\in S' (\rn)$ such that
\begin{\eq}   \label{2.9}
\| f \, | \Fs (\rn) \|_{\vp} = \Big\| \Big( \sum^\infty_{j=0} 2^{jsq} \big| (\vp_j \widehat{f})^\vee (\cdot) \big|^q \Big)^{1/q} | L_p (\rn) \Big\| < \infty
\end{\eq}
$($with the usual modification if $q= \infty)$.
\end{definition}

\begin{remark}   \label{R2.2}
As usual $\As (\rn)$ with $A \in \{ B,F \}$ means $\Bs (\rn)$ and $\Fs (\rn)$. 
The theory of these {\em inhomogeneous isotropic
spaces} and their history may be found in \cite{T83, T92, T06, T20}. In particular these spaces are independent of admitted resolutions
of unity $\vp$ according to \eqref{2.3}--\eqref{2.5} (equivalent quasi-norms). This justifies our omission of the subscript $\vp$ in 
\eqref{2.7}, \eqref{2.9} in the sequel. They are quasi--Banach spaces (Banach spaces if $p\ge 1$, $q\ge 1$).
\end{remark}
We recall some well--known properties and notation. Let $1<p<\infty$ and $k \in \no$. Then
\begin{\eq}   \label{2.10}
W^k_p (\rn) = F^k_{p,2} (\rn)
\end{\eq}
are the classical Sobolev spaces equivalently normed by
\begin{\eq}   \label{2.11}
\| f \, | W^k_p (\rn) \| = \Big( \sum_{|\alpha| \le k} \| D^\alpha f \, | L_p (\rn) \|^p \Big)^{1/p}.
\end{\eq}
including
\begin{\eq}   \label{2.12}
L_p (\rn) = F^0_{p,2} (\rn), \qquad 1<p<\infty.
\end{\eq}

Let
\begin{\eq}    \label{2.13}
w_\alpha (x) = (1+ |x|^2)^{\alpha/2}, \qquad \alpha \in \real, \quad x\in \rn.
\end{\eq}
Then $I_\alpha$,
\begin{\eq}   \label{2.14}
I_\alpha: \quad f \mapsto (w_\alpha \wh{f} )^\vee = (w_\alpha f^\vee )^\wedge, \qquad f\in S'(\rn), \quad \alpha \in \real.
\end{\eq}
is a lift  in the spaces $\As (\rn)$. This means,
\begin{\eq}   \label{2.15}
\begin{aligned}
I_\alpha \As (\rn) & = A^{s-\alpha}_{p,q} (\rn), \\
\big\| (w_\alpha \wh{f} )^\vee \, | A^{s-\alpha}_{p,q} (\rn) \big\| &\sim \|f \, |\As (\rn) \|.
\end{aligned}
\end{\eq}
One has in particular
\begin{\eq}   \label{2.16}
I_{-k} L_p (\rn) = W^k_p (\rn), \quad k\in \no, \quad 1<p<\infty,
\end{\eq}
with $W^k_p (\rn)$ as in \eqref{2.10}, \eqref{2.11}. The mappings in \eqref{1.5}, \eqref{1.6} are based on the special Besov spaces
\begin{\eq}   \label{2.17}
B^s_p (\rn) = B^s_{p,p} (\rn), \qquad s\in \real \quad \text{and} \quad 1 \le p \le \infty,
\end{\eq}
including the related H\"{o}lder--Zygmund spaces
\begin{\eq}   \label{2.18}
\Cc^s (\rn) = B^s_\infty (\rn), \qquad s \in \real.
\end{\eq}

Next we describe the dominating mixed counterpart of Definition \ref{D2.1} and Remark \ref{R2.2}. 
Let $\vp \in S(\real)$ with
\begin{\eq}   \label{2.19}
\vp (y) =1 \ \text{if $|y| \le 1$} \quad \text{and} \quad \vp(y) =0 \ \text{if $|y| \ge 3/2$}
\end{\eq}
and
\begin{\eq}   \label{2.20}
\vp_l (y) = \vp \big( 2^{-l} y \big) - \vp \big( 2^{-l+1} y \big), \qquad y \in \real, \quad l\in \nat,
\end{\eq}
be the one--dimensional version of \eqref{2.3}, \eqref{2.4}, notationally complemented by $\vp_0 = \vp$. Let
\begin{\eq}  \label{2.21}
\vp_k (x) = \prod^n_{l=1} \vp_{k_l} (x_l), \qquad k= (k_1, \ldots, k_n) \in \nat^n_0, \quad x = (x_1, \ldots, x_n )\in \rn.
\end{\eq}
Then
\begin{\eq}   \label{2.22}
\sum_{k \in \nat^n_0} \vp_k (x) =1, \qquad \text{for} \quad x\in \rn
\end{\eq}
is the counterpart of \eqref{2.5}. Furthermore $\vp_k$ refers to
\begin{\eq}    \label{2.23}
R_k = \big\{ x\in \rn: \ 2^{k_l -1} < |x_l| < 2^{k_l}, \ l=1, \ldots, n \big\}
\end{\eq}
($2^n$ rectangles). Recall that the entire analytic functions  $(\vp_k \wh{f} \, )^\vee (x)$ make sense pointwise in $\rn$ for
any $f\in S' (\rn)$. Let $[k] = \sum^n_{j=1}k_j$ if $k\in \zn$.

\begin{definition}  \label{D2.3}
Let $\vp$ and $\{ \vp_k \}_{k \in \nat^n_0}$ be the above dyadic resolution of unity.
\\[0.1cm]
{\em (i)} Let
\begin{\eq}   \label{2.24}
0<p \le \infty, \qquad 0<q \le \infty, \qquad r \in \real.
\end{\eq}
Then $\SB (\rn)$ is the collection of all $f\in S'(\rn)$ such that
\begin{\eq}   \label{2.25}
\| f \, | \SB (\rn) \|_{\vp} = \Big( \sum_{k\in \nat^n_0} 2^{r [k]q} \| (\vp_k \wh{f} )^\vee \, | L_p (\rn) \|^q \Big)^{1/q} <\infty
\end{\eq}
$($with the usual modification if $q=\infty)$.
\cm
{\em (ii)} Let
\begin{\eq}   \label{2.26}
0<p<\infty, \qquad 0<q \le \infty, \qquad r\in \real.
\end{\eq}
Then $\SF (\rn)$ is the collection of all $f \in S'(\rn)$ such that
\begin{\eq}   \label{2.27}
\| f \, | \SF (\rn) \|_{\vp} = \Big\| \Big( \sum_{k\in \nat^n_0} 2^{r[k]q} \, |(\vp_k \wh{f})^\vee (\cdot) |^q \Big)^{1/q} | L_p (\rn) \Big\| <\infty
\end{\eq}
$($with the usual modification if $q=\infty)$. 
\end{definition}

\begin{remark}   \label{R2.4}
As usual, $\SA (\rn)$ with $A \in \{B,F \}$ means $\SB (\rn)$ and $\SF (\rn)$. This is the dominating mixed counterpart of Definition
\ref{D2.1}. It coincides with \cite[Definition 1.4, p.~5]{T19}. In \cite[Remark 1.5, p.~5]{T19} one finds related references to the
books and papers already mentioned above, including in particular \cite{Schm80, Schm82, ST87, T10}. These spaces are independent of the
resolution of unity based on $\vp$ according to \eqref{2.19}--\eqref{2.22} (equivalent quasi--norms). This justifies our omission of the
subscript $\vp$ in \eqref{2.25} and \eqref{2.27} in the sequel. They are quasi--Banach spaces (Banach spaces if $p \ge 1$, $q \ge 1$).
Similarly as for the isotropic spaces in Remark \ref{R2.2} we recall some well--known properties and notation following \cite[Remark
1.5, p.~6]{T19} which in turn is based on related assertions in \cite{ST87}. Let $1<p<\infty$ and $r\in \no$. Then
\begin{\eq}   \label{2.28}
S^r_p W (\rn) = S^r_{p,2} F(\rn)
\end{\eq}
are the classical Sobolev spaces with dominating mixed smoothness which can be equivalently normed as in \eqref{1.3}, in particular,
\begin{\eq}   \label{2.29}
\| f \, | S^r_p W (\rn) \| = \sum_{\substack{\alpha \in \nat^n_0, \\ 0\le \alpha_j \le r}} \| D^\alpha f \, | L_p (\rn) \|
\end{\eq}
including
\begin{\eq}   \label{2.30}
L_p (\rn) = S^0_{p,2} F (\rn).
\end{\eq}
\end{remark}

Let $n\in \nat$ and
\begin{\eq}   \label{2.31}
v_\sigma (x) = \prod^n_{j=1} (1+ |x_j|^2 )^{\sigma/2}, \qquad \sigma \in \real, \quad x=(x_1, \ldots, x_n) \in \rn
\end{\eq}
be the counterpart of \eqref{2.13}. Then
\begin{\eq}  \label{2.32}
J_\sigma: \quad  f \mapsto (v_\sigma \wh{f} )^\vee = (v_\sigma f^\vee )^\wedge, \qquad f \in S'(\rn), \quad \sigma \in \real,
\end{\eq}
is a lift in the spaces  $\SA (\rn)$. This means for all $r\in \real$ and $0<p,q \le \infty$ ($p<\infty$ for $F$--spaces)
\begin{\eq}  \label{2.33}
\begin{aligned}
J_\sigma \SA (\rn) & = S^{r-\sigma}_{p,q}A (\rn), \\
\big\| (v_\sigma \wh{f} )^\vee \, | S^{r-\sigma}_{p,q}A (\rn) \big\| &\sim \|f \, |\SA (\rn) \|.
\end{aligned}
\end{\eq}
This is the dominating mixed counterpart of \eqref{2.15}. One has, similarly as in \eqref{2.16},
\begin{\eq}   \label{2.34}
J_{-r} L_p (\rn) = S^r_p W (\rn), \qquad r\in \no, \quad 1<p<\infty,
\end{\eq}
with $S^r_p W (\rn)$ as in  \eqref{2.29}, \eqref{2.30}. The mapping properties in the above Theorem are based on the special Besov
spaces with dominating mixed smoothness
\begin{\eq}   \label{2.35}
S^r_p B (\rn) = S^r_{p,p} B (\rn), \qquad r \in \real, \quad \text{and} \quad 1 \le p \le \infty,
\end{\eq}
including the related H\"{o}lder--Zygmund spaces
\begin{\eq}   \label{2.36}
S^r \Cc (\rn) = S^r_\infty B(\rn), \qquad r \in \real.
\end{\eq}
We refer the reader to \cite[pp.~5--7]{T19} for further equivalent quasi--norms and related references.

We rely in what follows mainly on wavelet representations and duality. Whereas  duality assertions will be described later on in the
course of the arguments we give now a detailed description of wavelet expansions for all spaces $\SA (\rn)$ according to Definition
\ref{D2.3}. We adapt corresponding assertions in \cite[pp.~14--16]{T19}, based in turn on \cite{Vyb06} and \cite{T10} to our later needs.

As usual, $C^{u} 
(\real)$ with $u \in \nat$ collects all complex-valued continuous functions on $\real$ having continuous bounded derivatives up to order $u$ inclusively. Let
\begin{\eq}   \label{2.37}
\psi_F \in C^{u} (\real), \quad \psi_M \in C^{u} (\real), \qquad u \in \nat,
\end{\eq}
be  real compactly supported Daubechies wavelets with
\begin{\eq}    \label{2.38}
\int_{\real} \psi_M (x) \, x^v \, \di x =0 \qquad \text{for all $v\in \no$ with $v<u$}
\end{\eq}
having $L_2$--norms 1. Let
\begin{\eq}   \label{2.39}
\psi_{-1,m} (t) = \sqrt{2} \, \psi_F (t-m) \quad \text{and} \quad \psi_{k,m} (t) = \psi_M (2^k t -m)
\end{\eq}
with $t\in \real$,  $k\in \no$ and $m\in \ganz$. Let $\nat_{-1} = \no \cup \{-1 \}$,
\begin{\eq}   \label{2.40}
\nat^n_{-1} = \{ k = (k_1, \ldots, k_n ), \ k_j \in \nat_{-1} \}, \qquad n \in \nat,
\end{\eq}
and
\begin{\eq}   \label{2.41}
\psi_{k,m} (x) = \prod^n_{j=1} \psi_{k_j, m_j} (x_j), \qquad k \in \nat^n_{-1}, \quad m \in \zn,
\end{\eq}
$x= (x_1, \ldots, x_n) \in \rn$. Then
\begin{\eq}   \label{2.42}
\big\{ 2^{[k]/2} \psi_{k,m}: \ k \in \nat^n_{-1}, \ m \in \zn \big\}
\end{\eq}
with (again) $[k] = \sum^n_{j=1} k_j$ is an orthonormal basis in $L_2 (\rn)$. Let $0< p,q \le \infty$ and $r\in \real$. Then $s^r_{p,q}
b(\rn)$ is the collection of all sequences
\begin{\eq}   \label{2.43}
\lambda = \big\{ \lambda_{k,m} \in \comp: \ k \in \nat^n_{-1}, \ m \in \zn \big\}
\end{\eq}
such that
\begin{\eq}  \label{2.44}
\| \lambda \, | s^r_{p,q} b (\rn) \| = \Big( \sum_{k \in \nat^n_{-1}} 2^{[k](r- \frac{1}{p})q}
\Big( \sum_{m \in \zn} |\lambda_{k,m} |^p \Big)^{q/p} \Big)^{1/q} < \infty
\end{\eq}
and $s^r_{p,q}f(\rn)$ is the collection of all sequences $\lambda$ in \eqref{2.43} such that
\begin{\eq}   \label{2.45}
\| \lambda \, |s^r_{p,q}f (\rn) \| = \Big\| \Big(\sum_{k,m} 2^{[k]rq}
|\lambda_{k,m} \, \chi_{k,m} (\cdot) |^q \Big)^{1/q} | L_p (\rn ) \Big\| < \infty
\end{\eq}
with the usual modifications if $p=\infty$ and/or $q=\infty$, where $\chi_{k,m}$ is the characteristic function of the rectangle
\begin{\eq}  \label{2.46}
R_{k,m} = \big\{ x=(x_1, \ldots, x_n) \in \rn, \ 2^{-k_j} m_j <x_j < 2^{-k_j} (m_j +1) \big\}.
\end{\eq}
Let 
\begin{\eq}   \label{2.47}
\sigma_p = \max \big( \frac{1}{p}, 1) -1 \quad \text{and} \quad \sigma_{p,q} = \max \big( \frac{1}{p}, \frac{1}{q}, 1 \big) -1.
\end{\eq}

\begin{proposition}    \label{P2.5}
{\em (i)} Let $0<p \le \infty$, $0<q \le \infty$ and $r\in \real$. Let $\psi_{k,m}$ be the wavelets in \eqref{2.41} based on \eqref{2.37},
\eqref{2.38} with $u> \max (r, \sigma_p -r)$. Let $f\in S'(\rn)$. Then $f\in \SB (\rn)$ if, and only if, it can be represented as
\begin{\eq}   \label{2.48}
f = \sum_{k \in \nat^n_{-1}} \sum_{m \in \zn} \lambda_{k,m} \, \psi_{k,m}, \qquad \lambda \in s^r_{p,q}b (\rn).
\end{\eq}
The representation \eqref{2.48} is unique,
\begin{\eq}   \label{2.49}
\lambda_{k,m} = \lambda_{k,m} (f) = 2^{[k]} (f, \psi_{k,m}), \qquad k \in \nat^n_{-1}, \quad m\in \zn,
\end{\eq}
and
\begin{\eq}   \label{2.50}
J: \quad f \mapsto \big\{ \lambda_{k,m}: \ k\in \nat^n_{-1}, \ m\in \zn \big\}
\end{\eq}     
is an isomorphic mapping of $\SB (\rn)$ onto $s^r_{p,q} b (\rn)$.
\cm
{\em (ii)} Let $0<p<\infty$, $0<q \le \infty$ and $r\in \real$. Let $u > \max (r, \sigma_{p,q} -r)$. Let $f\in S'(\rn)$. Then $f \in
\SF (\rn)$ if, and only if, it can be represented as
\begin{\eq}   \label{2.51}
f = \sum_{k \in \nat^n_{-1}} \sum_{m \in \zn} \lambda_{k,m} \, \psi_{k,m}, \qquad \lambda \in s^r_{p,q}f (\rn).
\end{\eq}
The representation \eqref{2.51} is unique with \eqref{2.49}. Furthermore, $J$ in \eqref{2.50} is an isomorphic mapping of $\SF (\rn)$
onto $s^r_{p,q} f (\rn)$.
\end{proposition}

\begin{remark}   \label{R2.6}
This is a modified version of \cite[Theorem 1.12, pp.~15, 16]{T19} $($with different normalizations$)$. There one finds further 
explanations and references, especially to \cite{Vyb06, T10}. Of special interest for us will be the case
\begin{\eq}   \label{2.52}
S^r \Cc (\rn) = S^r_\infty B (\rn) = S^r_{\infty, \infty} B (\rn), \qquad r\in \real.
\end{\eq}
Then one has the representation
\begin{\eq}   \label{2.53}
f = \sum_{k\in \nat^n_{-1}, m\in \zn} \lambda_{k,m} \, \psi_{k,m}, \qquad \lambda \in s^r_{\infty,\infty} b (\rn)
\end{\eq}
and
\begin{\eq}   \label{2.54}
\begin{aligned}
\|f \, | S^r \Cc (\rn) \| & \sim \| \lambda \, | s^r_{\infty,\infty} b (\rn) \| \\
& = \sup_{k\in \nat^n_{-1}, m\in \zn} 2^{[k]r} \, |\lambda_{k,m}|, \qquad r \in \real.
\end{aligned}
\end{\eq}
\end{remark}

\section{Proof of the Theorem}   \label{S3}
We break the proof of the Theorem in 10 steps.
\cm
{\em Step 1.} First we prove
\begin{\eq}   \label{3.1}
F: \quad S^0 \Cc (\rn) \hra S^{-1} \Cc (\rn)
\end{\eq}
for the H\"{o}lder--Zygmund spaces according to \eqref{2.36}. Let $f\in S^0 \Cc (\rn)$. We expand $Ff \in S' (\rn)$ according to
Proposition \ref{P2.5},
\begin{\eq}   \label{3.2}
Ff = \sum_{k\in \nat^n_{-1}} \sum_{m\in \zn} \lambda_{k,m} \, \psi_{k,m}
\end{\eq}
with
\begin{\eq}    \label{3.3}
\lambda_{k,m} = 2^{[k]} (Ff, \psi_{k,m}) = 2^{[k]} (f, F \psi_{k,m} )
\end{\eq}
where we used that $F$ is self--dual in the context of the framework of the dual pairing $\big( S(\rn), S'(\rn) \big)$, $F' = F$. This
possibility is not totally obvious, but we add a comment about this representability of $Ff$ in Remark \ref{R3.1} below. Let 
temporarily $F_1$ be the Fourier transform on the real line \real. Then one obtains  from \eqref{2.41} that
\begin{\eq}   \label{3.4}
(F \psi_{k,m} )(x) = \prod^n_{j=1}  (F_1 \psi_{k_j, m_j} )(x_j), \qquad k \in \nat^n_{-1}, \quad m\in \zn,
\end{\eq}
$x = (x_1, \ldots, x_n )$. One has by \cite[Proposition 1.19, p.~19]{T19}
\begin{\eq}   \label{3.5}
\| F \psi_{k,m} \, | S^0_1 B(\rn) \| = \prod^n_{j=1} \| F_1 \psi_{k_j, m_j} \, |B^0_1 (\real \|
\end{\eq}
and by \cite[Step1 of the proof of Theorem 2.12]{HST24} (or by direct calculation) that
\begin{\eq}   \label{3.6}
\| F_1 \psi_{k_j, m_j} \, |B^0_1 (\real) \| \le c, \qquad k_j \in \nat_{-1}, \quad m_j \in \ganz,
\end{\eq}
for some $c >0$ (independently of $k_j$ and $m_j$). Using in addition the duality
\begin{\eq}    \label{3.7}
S^0_1 B(\rn)' = S^0_\infty B(\rn) = S^0 \Cc (\rn)
\end{\eq}
in \cite[Proposition 1.17, pp.~17--18]{T19} it follows from \eqref{3.3}, \eqref{3.5}, \eqref{3.6} that
\begin{\eq}  \label{3.8}
\begin{aligned}
|\lambda_{k,m} | & \le 2^{[k]} \|f \, |S^0 \Cc (\rn) \| \cdot \|F \psi_{k,m} \, | S^0_1 B(\rn) \| \\
&\le c \, 2^{[k]} \| f \, |S^0 \Cc (\rn) \|.
\end{aligned}
\end{\eq}
Now one obtains \eqref{3.1} from \eqref{2.52}--\eqref{2.54} with $r=-1$.
\cm
{\em Step 2.} Let $2 \le p \le \infty$. We justify
\begin{\eq}   \label{3.9}
F: \quad S^0_p B(\rn) \hra S^{\frac{2}{p} -1}_p B (\rn), \qquad 2 \le p \le \infty,
\end{\eq}
by complex interpolation $[\cdot, \cdot]_\theta$, $0<\theta <1$, of \eqref{3.1},
\begin{\eq}   \label{3.10}
F: \quad S^0_\infty B(\rn) \hra S^{-1}_\infty B(\rn),
\end{\eq}
and
\begin{\eq}   \label{3.11}
F: \quad L_2 (\rn) \hra L_2 (\rn),
\end{\eq}
rewritten according to \eqref{2.30} and $S^r_{p,p} B (\rn) = S^r_{p,p} F(\rn)$, $1<p<\infty$, as
\begin{\eq}   \label{3.12}
F: \quad S^0_2 B(\rn) \hra S^0_2 B(\rn).
\end{\eq}
Based on the isomorphic mappings of these spaces onto corresponding sequence spaces as described in Proposition \ref{P2.5} one can 
shift this task to the complex interpolation of the related sequence spaces. But then one is  essentially in the same situation as
for isotropic spaces  in \cite{HST23} and \cite[Step 2 of the proof of Theorem 2.12]{HST24} with the outcome
\begin{\eq}   \label{3.13}
[S^0\Cc (\rn), L_2 (\rn) ]_\theta = [S^0_\infty B(\rn), S^0_2 B (\rn) ]_\theta = S^0_p B (\rn)
\end{\eq}
and
\begin{\eq}   \label{3.14}
[S^{-1} \Cc (\rn, L_2 (\rn) ]_\theta = [S^{-1}_\infty B (\rn), S^0_2 B(\rn) ]_\theta = S^{\frac{2}{p} -1}_p B (\rn)
\end{\eq}
where $0<\theta <1$,
\begin{\eq}   \label{3.15}
\frac{1}{p} = \frac{1-\theta}{\infty} + \frac{\theta}{2} = \frac{\theta}{2}
\quad \text{and} \quad \theta -1 = \frac{2}{p} -1 = \frac{1}{p} - \frac{1}{p'}
\end{\eq}
with $\frac{1}{p} + \frac{1}{p'} =1$. Now \eqref{3.10}--\eqref{3.14} prove \eqref{3.9}.
\cm
{\em Step 3.} Let $1 \le p \le2$. We justify
\begin{\eq}    \label{3.16}
F: \quad S^{\frac{2}{p} -1}_p B(\rn) \hra S^0_p B (\rn), \qquad 1 \le p \le 2,
\end{\eq}
by duality. By \cite[Proposition 1.17, pp.~17--18]{T19} one has the duality
\begin{\eq}   \label{3.17}
S^r_p B (\rn)' = S^{-r}_{p'} B(\rn), \qquad 1 \le p <\infty, \quad \frac{1}{p} + \frac{1}{p'} =1, \quad r\in \real,
\end{\eq}
in the framework of the dual pairing $\big( S(\rn), S'(\rn) \big)$ complemented by
\begin{\eq}   \label{3.18}
\Si{}^r_\infty B (\rn)' = S^{-r}_1 B(\rn), \qquad r\in \real.
\end{\eq}
Here $\Si{}^r_\infty B (\rn) = \Si{}^r \Cc (\rn)$ is the completion  of $S(\rn)$ in $S^r \Cc (\rn)$. Recall that $F$ is self--dual, 
$F = F'$. Then it follows by duality from \eqref{3.9}, $\frac{2}{p} -1 = \frac{1}{p}- \frac{1}{p'}$  and \eqref{3.10}, complemented by
\begin{\eq}   \label{3.19}
F: \quad \Si{}^0_\infty B (\rn) \hra \Si{}^{-1}_\infty B(\rn)
\end{\eq}
that
\begin{\eq}  \label{3.20}
F: \quad S^{\frac{1}{p'} - \frac{1}{p}}_{p'} B (\rn) \hra S^0_{p'} B (\rn), \qquad 2 \le p \le \infty, \quad \frac{1}{p}- \frac{1}{p'}
=1.
\end{\eq}
This proves \eqref{3.16}.
\cm
{\em Step 4.} The continuity of the mapping \eqref{1.7}, \eqref{1.8} follows now from \eqref{3.16}, \eqref{3.9} and the monotonicity
of the spaces $S^r_p B (\rn)$ with respect to $r$ for fixed $p$.
\cm
{\em Step 5.} We justify in three steps that  the mapping in \eqref{1.7}, \eqref{1.8} with $r_1 \ge 0$ and $r_2 \ge 0$ are not compact
if either $r_1 =0$ or $r_2 = 0$. This requires some efforts. First we show that
\begin{\eq}  \label{3.21}
F: \quad S^{r_1}_p B (\rn) \hra S^{\frac{2}{p} -1}_p B (\rn), \qquad 2 \le p \le \infty, \quad r_1 \ge 0,
\end{\eq}
is not compact. By $\frac{2}{p} - 1 - \frac{1}{p} = - \frac{1}{p'}$ and the embedding
\begin{\eq}   \label{3.22}
\id: \quad S^{\frac{2}{p} -1}_p B (\rn) \hra S^{- \frac{1}{p'}} \Cc (\rn) = S^{- \frac{1}{p'}}_\infty B (\rn)
\end{\eq}
according to \cite[(1.312), p.~55]{T19} it is sufficient to prove that
\begin{\eq}   \label{3.23}
F: \quad S^{r_1}_p B(\rn) \hra S^{- \frac{1}{p'}} \Cc (\rn), \qquad 2 \le p \le \infty, \quad r_1 \ge 0,
\end{\eq}
is not compact. Let $\psi$ be a non--trivial $C^\infty$ function in \real,
\begin{\eq}  \label{3.24}
\psi_k (x) = \prod^n_{j=1} \psi (2^{-k_j} x_j), \qquad k= (k_1, \ldots, k_n ) \in \nat^n_0, \quad x=(x_1, \ldots, x_n) \in \rn,
\end{\eq}
such that
\begin{\eq}   \label{3.25}
\psi_k (x) \vp_k (x) = \psi_k (x) \qquad \text{if} \quad k= (k_1, 0, \ldots,0), \quad k_1 \in \no,
\end{\eq} 
and
\begin{\eq}   \label{3.26}
\supp \psi_k \cap \supp \vp_l = \emptyset \qquad \text{if} \quad l \not= k,
\end{\eq}
where $\vp_k$ and $\vp_l$ are the same functions as in the Definition \ref{D2.3}. Let
\begin{\eq}   \label{3.27}
f_k (x) = 2^{-\frac{k_1}{p}} \psi_k (x), \qquad \text{$k$ as above}, \quad x\in \rn.
\end{\eq}
Then one has for the Sobolev spaces in \eqref{2.28}, \eqref{2.29} (notationally extended to $p=\infty$)
\begin{\eq}   \label{3.28}
\| f_k \, | S^r_p W (\rn) \| = \sum_{\substack{\alpha \in \nat^n_0, \\ 0\le \alpha_j \le r}} \| D^\alpha f_k \, | L_p (\rn) \| 
\sim 1, \qquad k_1 \in \nat,
\end{\eq}
$r\in \nat$, $1<p \le \infty$. Using the monotonicity as described in \cite[Section 1.3.2, pp.~55--56]{T19}, especially
\begin{\eq}   \label{3.29}
S^r_p W (\rn) \hra S^{r_1}_p B (\rn) \hra L_p (\rn), \qquad 0<r_1 <r, \quad 2 \le p \le \infty,
\end{\eq}
it follows
\begin{\eq}   \label{3.30}
\|f_k \, | S^{r_1}_p B (\rn) \| \sim 1, \qquad 2 \le p \le \infty, \quad r_1 \ge 0.
\end{\eq}
One has by Definition \ref{D2.3}
\begin{\eq}   \label{3.31}
\begin{aligned}
\| \wh{f_k} \, | S^{- \frac{1}{p'}} \Cc (\rn) \| &\sim \sup_{l\in \nat^n_0, x\in \rn} 2^{-\frac{k_1}{p'}} \big| (\vp_l f_k)^\vee (x) \big| \\
&\sim \sup_{x\in \rn} 2^{- \frac{k_1}{p'} - \frac{k_1}{p}} \big| \psi^\vee_k (x) \big| \sim 1
\end{aligned}
\end{\eq}
and
\begin{\eq}   \label{3.32}
\| \wh{f_{k^1}} - \wh{f_{k^2}} \, | S^{-\frac{1}{p'}} \Cc (\rn) \| \sim 1, \qquad k^j = (k^j_1, 0, \ldots, 0),
\end{\eq}
$k^1 \not= k^2$. This shows that $F$ in \eqref{3.23} is not compact.
\cm
{\em Step 6}. We prove that
\begin{\eq}   \label{3.33}
F: \quad S^0_p B(\rn) \hra S^{\frac{2}{p} -1 - r_2}_p B (\rn), \qquad 2 \le p \le \infty, \quad r_2 \ge 0,
\end{\eq}
is not compact. It follows from an embedding as in \eqref{3.22} that it is sufficient to show  that
\begin{\eq}   \label{3.34}
F: \quad  S^0_p B (\rn) \hra S^\sigma \Cc (\rn) = S^\sigma_\infty B (\rn), \qquad 2 \le p \le \infty, \quad \sigma \le -1,
\end{\eq}
is not compact. One can take over the related arguments from \cite[proof of Theorem 2.12]{HST24}, based on \cite{HST23}, for isotropic
spaces. Let $\psi_F (t)$ be a smooth compactly supported father wavelet according to \eqref{2.37}. Let $\psi(x) = \prod^n_{j=1}
\psi_F (x_j)$ which refers to $k= (-1, \ldots, -1)$ in \eqref{2.40}. Then $\psi_m (x) = \psi(x-m)$, $x\in \rn$, $m\in \zn$ are wavelets for the isotropic spaces and also for the above spaces with dominating mixed smoothness. Let
\begin{\eq}   \label{3.35}
f_m (x) = (F^{-1} \psi_m )(x) = e^{imx} (F^{-1} \psi)(x), \qquad x\in \rn, \quad m\in \zn.
\end{\eq}
Then $\{ Ff_m = \psi_m: \ m\in \zn \}$ is not compact in any space $S^\sigma \Cc (\rn)$, $\sigma \in \real$. We wish to show that
uniformly 
\begin{\eq}   \label{3.36}
\| f_m \, |S^0_p B (\rn) \| \sim 1, \qquad m\in \zn,
\end{\eq}
based on \eqref{2.25}. There are two cases. One has either
\begin{\eq}   \label{3.37}
(\vp_k \wh{f_m})^\vee (x) = (\vp_k \psi_m )^\vee (x) =0
\end{\eq}
or (assuming $\vp_k \psi_m = \psi_m$)
\begin{\eq}   \label{3.38}
\| (\vp_k \wh{f_m} )^\vee (\cdot) \, |L_p (\rn) \| = \|e^{imx} \psi^\vee (\cdot) \, |L_p (\rn) \| = \| \psi^\vee \, | L_p(\rn) \|
\end{\eq}
for one $k$. Inserted in \eqref{2.25} one obtains \eqref{3.36}. This shows that $F$ in \eqref{3.34} and \eqref{3.33} is not compact.
\cm
{\em Step 7.} The two preceding steps show that $F$ in \eqref{1.8},
\begin{\eq}   \label{3.39}
F: \quad S^{r_1}_p B (\rn) \hra S^{\frac{1}{p} - \frac{1}{p'} -r_2}_p B(\rn), \qquad 2 \le p \le \infty,
\end{\eq}
is not compact if either $r_1 = 0$ or $r_2 =0$. The corresponding assertion for $1 \le p \le 2$ can be obtained by duality as follows.
Let $F$ in \eqref{1.7},
\begin{\eq}   \label{3.40}
F: \quad S^{\frac{1}{p} - \frac{1}{p'} + r_1}_p B (\rn) \hra S^{-r_2}_p B (\rn), \qquad 1\le p \le 2,
\end{\eq}
$r_1 \ge 0$, $r_2 \ge 0$, be compact.  Then it follows from \eqref{3.17} and $F =F'$ that
\begin{\eq}   \label{3.41}
F: \quad S^{r_2}_{p'} B (\rn) \hra S^{\frac{1}{p'} - \frac{1}{p} - r_1}_{p'} B (\rn)
\end{\eq}
is compact. But this is not the case if either $r_1 =0$ or $r_2 =0$.
\cm
{\em Step 8}. In the preceding steps we adapted the corresponding proofs in \cite[Theorem 2.12]{HST24} and the underlying paper 
\cite{HST23} for related isotropic spaces to their dominating mixed counterparts. This is not possible any longer in order to justify
that the mapping $F$ in \eqref{1.7}, \eqref{1.8} is compact if both $r_1 >0$ and $r_2 >0$. In \cite{HST23,HST24} we relied on an
elaborated theory for weighted isotropic spaces $\As (\rn, w)$ with suitable weights $w$, including compact mapping properties 
expressed in terms of entropy numbers, and corresponding wavelet representations. There is (as far as we know) no counterpart for
related weighted spaces $\SA (\rn,w)$. We circumvent this shortcoming and shift the problem of compactness (not of entropy numbers) 
from spaces with dominating mixed smoothness to appropriate isotropic spaces. This requires again some efforts and will be done in the
present step, preparing the proof of compactness in the next step. First we follow \cite[Section 6.1, pp.~263--268]{T06}. The class
$W^n$, $n\in \nat$, of admissible weight functions is the collection of all positive $C^\infty$ functions in $\rn$ such that for all
$\gamma \in \nat^n_0$ and some $c_\gamma >0$,
\begin{\eq}  \label{3.42}
|D^\gamma w(x)| \le c_\gamma \, w(x), \qquad x\in \rn,
\end{\eq}
and for some $\beta \ge 0$ and $c>0$,
\begin{\eq}   \label{3.43}
0< w(x) \le c \, w(y) \big( 1 + |x-y|^2 \big)^{\beta/2}, \qquad x\in \rn, \quad y \in \rn.
\end{\eq}
Both weights
\begin{\eq}  \label{3.44}
w_\alpha (x) = (1 + |x|^2 )^{\alpha/2}, \qquad \alpha \in \real, \quad x \in \rn,
\end{\eq}
according to \eqref{2.13} and
\begin{\eq}   \label{3.45}
v_\sigma (x) = \prod^n_{j=1} (1 + |x_j|^2)^{\sigma/2}, \qquad \sigma \in \real, \quad x = (x_1, \ldots, x_n) \in \rn,
\end{\eq}
according to \eqref{2.31} (as the product of related one--dimensional weights) belong to $W^n$. For these weights  we introduced in
\cite{T06} (based on the references given there) the weighted spaces $\As (\rn, w)$ replacing in Definition \ref{D2.1} the Lebesgue
spaces $L_p (\rn)$ by $L_p (\rn,w)$ with
\begin{\eq}   \label{3.46}
\| f \, | L_p (\rn, w)\| = \| wf \, | L_p (\rn) \|, \qquad 0<p \le \infty, \quad w\in W^n.
\end{\eq}
It is one of the main observations that $f \mapsto wf$ is an isomorphic mapping of $\As (\rn, w)$ onto $\As (\rn)$,
\begin{\eq}    \label{3.47}
\| wf \, | \As (\rn) \| \sim \|f \, | \As (\rn,w) \|, \quad A\in \{B,F\}, \quad s\in \real, \ 0<p,q \le \infty
\end{\eq}
($p<\infty$ for $F$--spaces). Details and references may be found in \cite[Section 6.1.3, pp.~265--266]{T06}. In particular, for fixed
$p$ and $q$ the spaces $\As (\rn, w)$ are monotonically included both for the smoothness $s$ and for related weights. There is little
doubt that there are suitable counterparts for the spaces $\SA (\rn)$ as introduced in Definition \ref{D2.3}, at least for weights
of type \eqref{3.45}. But this has not yet been done. It will be sufficient for us to introduce the spaces $\SA (\rn, v_\sigma)$ with
$v_\sigma$ in \eqref{3.45} as the collection of all $f\in S'(\rn)$ such that $v_\sigma f \in \SA (\rn)$,
\begin{\eq}   \label{3.48}
\| f \, | \SA (\rn, v_\sigma) \| = \| v_\sigma f \, | \SA (\rn) \|.
\end{\eq}
For fixed $p$ and $q$ there is still the monotonicity  with respect to the smoothness $r$. But we shift the missing monotonicity with
respect to the weights to the isotropic case. This is based on the following observations. Spaces with dominating mixed smoothness
have been compared with other spaces in $\rn$ in detail. We refer the reader to \cite{Schm07} and \cite{NgS17a, NgS17b}. For our 
purpose the following specific assertions are sufficient. Let $B^s_p (\rn)$, $\Cc^s (\rn) = B^s_\infty (\rn)$ be as in \eqref{2.17},
\eqref{2.18} and $S^r_p B(\rn)$, $S^r \Cc (\rn) = S^r_\infty B (\rn)$ be as in \eqref{2.35},
\eqref{2.36}. Let $2 \le n \in \nat$. Then
\begin{\eq}  \label{3.49}
B^{rn}_p (\rn) \hra S^r_p B (\rn) \hra B^r_p (\rn), \qquad 1 \le p \le \infty, \quad r>0,
\end{\eq}
and
\begin{\eq}   \label{3.50}
S^0_p B(\rn) \hra B^0_p (\rn), \qquad 1 \le p \le 2.
\end{\eq}
These are special cases of  \cite[Theorems 3.1, 3.6]{NgS17a}. 
\cm
{\em Step 9.} After the preparations in Step 8 we prove now that the mappings $F$ in \eqref{1.7}, \eqref{1.8} are compact if both
$r_1 >0$ and $r_2 >0$. Let $1 \le p \le 2$ and $r_1 >0$. Then it follows from the lifting property \eqref{2.33} that
\begin{\eq}   \label{3.51}
\| (v_{r_1} \wh{f} )^\vee | S^{\frac{1}{p} - \frac{1}{p'}}_p B (\rn) \| \sim \| f \, | S^{\frac{1}{p} - \frac{1}{p'} +r_1}_p B(\rn) \|
\end{\eq}
and by Step 3 that
\begin{\eq}   \label{3.52}
\| v_{r_1} \wh{f} \, |S^0_p B(\rn) \| \le c \, \|f \, |S^{\frac{1}{p} - \frac{1}{p'} +r_1}_p B (\rn) \|.
\end{\eq}
This means in terms of the preceding Step 8, that both
\begin{\eq}   \label{3.53}
F: \quad S^{\frac{1}{p} - \frac{1}{p'} +r_1}_p B(\rn) \hra S^0_p B ( \rn, v_{r_1}), \qquad 1\le p \le 2, \quad r_1 >0,
\end{\eq}
and, using in addition \eqref{3.50},
\begin{\eq}   \label{3.54}
F: \quad S^{\frac{1}{p} - \frac{1}{p'} +r_1}_p B (\rn) \hra B^0_p (\rn, v_{r_1} ), \qquad 1 \le p \le 2, \quad r_1 > 0,
\end{\eq}
are continuous. It follows from \eqref{3.44}, \eqref{3.45} that
\begin{\eq}   \label{3.55}
w_\alpha (x) \le n^{\alpha/2} v_\alpha (x) \quad \text{and} \quad v_\alpha (x) \le w_{n\alpha} (x), \qquad \alpha >0, \quad x\in \rn.
\end{\eq}
As mentioned in Step 8 the spaces $\Bs (\rn,w)$ are monotonically included with respect to the weights (for fixed $s,p,q$). Then one 
has by \eqref{3.54} that  also the mapping
\begin{\eq}   \label{3.56}
F: \quad S^{\frac{1}{p} - \frac{1}{p'} + r_1}_p B(\rn) \hra B^0_p (\rn, w_{r_1} ), \qquad 1 \le p \le 2, \quad r_1 >0,
\end{\eq}
is continuous. Combined with the compact embedding
\begin{\eq}   \label{3.57}
\id: \quad B^0_p (\rn, w_{r_1} ) \hra  B^{-r_2}_p (\rn), \quad 1 \le p \le 2, \quad r_1 >0, \quad r_2 >0,
\end{\eq}
as a very special case of \cite[Theorem 6.31, pp.~282--283]{T06}, one obtains that
\begin{\eq}  \label{3.58}
F: \quad S^{\frac{1}{p} - \frac{1}{p'} +r_1}_p B (\rn) \hra B^{-r_2}_p (\rn), \quad 1 \le p \le 2, \quad r_1 >0, \quad r_2 >0,
\end{\eq}
is compact. By duality as already used in Step 3 (and its obvious counterpart for isotropic spaces) it follows from \eqref{3.58}
that also
\begin{\eq}   \label{3.59}
F: \quad B^{r_2}_p (\rn) \hra S^{\frac{1}{p} - \frac{1}{p'} -r_1}_p B(\rn), \quad 2 \le p \le \infty, \quad r_1 >0, \quad r_2 >0,
\end{\eq}
is compact.  Then \eqref{3.49} (and appropriate adaption of the parameters) show that
\begin{\eq}  \label{3.60}
F: \quad S^{r_1}_p B (\rn) \hra S^{\frac{1}{p} - \frac{1}{p'} -r_2}_p B(\rn), \quad 2 \le p \le \infty, \quad r_1 >0, \quad r_2 >0,
\end{\eq}
and
\begin{\eq}  \label{3.61}
F: \quad \os{\circ}{S}{}^{r_1} \Cc (\rn) \hra \os{\circ}{S}{}^{-1 - r_2} \Cc (\rn), \quad r_1 >0, \quad r_2 >0,
\end{\eq}
are compact. This proves that $F$ in \eqref{1.8} is compact if both $r_1 >0$ and $r_2 >0$. Application of duality again shows that also
the mapping in \eqref{1.7} is compact if both $r_1 >0$ and $r_2 >0$.
\cm
{\em Step 10.} It remains to prove that there are no continuous mappings as in \eqref{1.7}, \eqref{1.8} if either $r_1 <0$ or $r_2 <0$.
Let us assume that there is a continuous mapping
\begin{\eq}   \label{3.62}
F: \quad S^{r_1}_p B (\rn) \hra S^{\frac{2}{p} -1 -r_2}_p B(\rn) \quad \text{for some $r_1 <0$, $r_2 \in \real$, $2 \le p \le \infty.$}
\end{\eq}
By embedding it is sufficient to deal with the case $r_2 >0$. We already know that 
\begin{\eq}   \label{3.63}
F: \quad S^{r_0}_p B(\rn) \hra S^{\frac{2}{p} -1-r_2}_p B(\rn) \quad r_0 >0, \quad r_2 >0, \quad 2 \le p \le \infty,
\end{\eq}
is compact. We use the complex interpolation $[\cdot, \cdot]_\theta$, $0<\theta <1$, according to  \cite[Proposition 1.20, p.~20]{T19}
and the references given there. One has as a special case
\begin{\eq}   \label{3.64}
\big[ S^{r_0}_p B (\rn), S^{r_1}_p B (\rn)\big]_\theta = S^0_p B(\rn), \qquad 0= (1-\theta)r_0 + \theta r_1,
\end{\eq}
$2 \le p <\infty$. However if $F$ in \eqref{3.62} is continuous then the interpolation \eqref{3.64} with the compact mapping 
\eqref{3.63} shows that also
\begin{\eq}   \label{3.65}
F: \quad S^0_p B(\rn) \hra S^{\frac{2}{p}-1 -r_2}_p B (\rn)
\end{\eq}
is compact. This follows from the interpolation theory for compact operators  as described in \cite[Section 1.16.4, pp.~117--118]{T78}.
But we know by Step 6 that this is not the case. This disproves that there is a continuous mapping as assumed in \eqref{3.62}. 
Similarly one can argue in the other cases with $1\le p <\infty$. If $p=\infty$ then it follows from the assumption
\begin{\eq}  \label{3.66}
F: \quad S^{r_1} \Cc (\rn) \hra S^{-1 -r_2} \Cc (\rn), \qquad r_1 <0, \quad r_2 \in \real
\end{\eq} 
that also
\begin{\eq}    \label{3.67}
F: \quad \os{\circ}{S}{}^{r_1} \Cc (\rn) \hra \os{\circ}{S}{}^{-1-r_2} \Cc (\rn), \qquad r_1 <0, \quad r_2 \in \real.
\end{\eq}
The duality \eqref{3.18} shifts this question  to the case $p=1$ which is covered by the above arguments. This shows that \eqref{3.66}
requires $r_1 \ge 0$ and $r_2 \ge 0$.

\begin{remark}   \label{R3.1}
We add a comment in connection of the representability of $Ff$ in \eqref{3.2}, \eqref{3.3}. Let $B^s_p (\rn, w_\alpha) = B^s_{p,p}
(\rn, w_\alpha)$, $1\le p \le \infty$, $\alpha \in \real$, be the related isotropic spaces as considered in Step 8 of the above proof. 
Then one has for fixed $p$,
\begin{\eq}    \label{3.68}
S(\rn) = \bigcap_{\alpha \in \real, s\in \real} B^s_p (\rn, w_\alpha) \quad \text{and} \quad S'(\rn) = \bigcup_{\alpha \in \real, s\in \real}
B^s_p (\rn, w_\alpha).
\end{\eq}
This (more or less obvious) assertion may be found in \cite[Remark 2.91, p.~74]{T20} with a reference to  \cite{Kab08}. There is 
little doubt that there is a counterpart for the spaces $S^r_p B(\rn, v_\sigma) = S^r_{p,p} B(\rn, v_\sigma)$ as introduced in 
\eqref{3.48},
\begin{\eq}   \label{3.69}
S(\rn) = \bigcap_{r \in \real, \sigma\in \real} S^r_p B (\rn, v_\sigma) \quad \text{and} \quad S'(\rn) = \bigcup_{r \in \real, \sigma\in \real} S^r_p B (\rn, v_\sigma).
\end{\eq}
But this is not (yet) available in the literature. The proof of the isotropic version \cite[Theorem 2.12]{HST24}, based on 
\cite{HST23}, of the above Theorem relies on the possibility to expand any element of $B^s_p (\rn, w_\alpha)$ and, by \eqref{3.68}
of any element  of $S'(\rn)$ in terms of Daubechies wavelets. One can expect that there is a dominating mixed counterpart based on the
extension of Proposition \ref{P2.5} to the weighted spaces $\SA (\rn, v_\sigma)$. But again this is not (yet) available in the 
literature. However one can circumvent these technical shortcomings in the same way as in the above proof (avoiding weights).
Let $g$ be the right--hand side of \eqref{3.2} with \eqref{3.3}. Then it follows from \eqref{3.68} that both $g$ and $Ff$ (multiplied 
with suitable cut--off functions) belong locally to, say, some space $B^s_2(\rn)$. Let $h \in B^{-s}_2 (\rn) = B^s_2 (\rn)'$ with
compact support. Then it follows from Proposition \ref{P2.5} and
\begin{\eq}   \label{3.70}
(Ff, \psi_{k,m} ) = ( g, \psi_{k,m}), \qquad k \in \nat^n_{-1}, \quad m\in \zn,
\end{\eq}
that $(Ff,h) = (g,h)$, and in particular,
\begin{\eq}   \label{3.71}
(Ff, \vp ) = (g, \vp), \qquad \vp \in D(\rn) = C^\infty_0 (\rn).
\end{\eq}
But this ensures $Ff =g$ in the framework of the dual pairing $\big(S(\rn), S'(\rn) \big)$ and also the representability \eqref{3.2},
\eqref{3.3}
\end{remark}

\begin{problem}   \label{P3.2}
The above arguments show again that it would be desirable to develop a theory for the weighted  spaces $\SA (\rn, v_\sigma)$ in full
generality parallel at least to corresponding  assertions for the weighted isotropic spaces $\As (\rn, w_\alpha)$. Then one can
avoid the somewhat artistic arguments in the Steps 8 and 9 and in the above Remark \ref{R3.1}. In addition to assertions being similar
to corresponding ones for weighted isotropic spaces there might be also some new aspects, based on Faber bases, for numerical 
integration and discrepancy (number--theoretical properties  about the distribution of points).
\end{problem}

\end{document}